\documentclass[12pt,reqno]{article}

\usepackage[usenames]{color}
\usepackage{amsmath}
\usepackage{amsfonts}
\usepackage{amssymb}
\usepackage{graphicx}
\usepackage{amscd}
\usepackage{multirow}
\usepackage{color}
\usepackage{psfig}
\usepackage{graphics,amsmath,amssymb}
\usepackage{amsthm}
\usepackage{amsfonts}
\usepackage{latexsym}
\usepackage{epsf}

\graphicspath{%
    {converted_graphics/}
    {C:/}
}
\begin{document}

\theoremstyle{plain}
\newtheorem{theorem}{Theorem}
\newtheorem{corollary}[theorem]{Corollary}
\newtheorem{lemma}[theorem]{Lemma}
\newtheorem{proposition}[theorem]{Proposition}

\theoremstyle{definition}
\newtheorem{definition}[theorem]{Definition}
\newtheorem{example}[theorem]{Example}
\newtheorem{conjecture}[theorem]{Conjecture}

\theoremstyle{remark}
\newtheorem{remark}[theorem]{Remark}


\phantom{\vspace{1cm}}
\phantom{\vspace{1cm}}
\phantom{\vspace{1cm}}

\begin{center}
\vskip 1cm
{\LARGE\bf Maximal Gaps Between Prime $k$-Tuples: \\ A Statistical Approach}
\vskip 1cm
\large
Alexei Kourbatov\\
JavaScripter.net/math\\
{\tt akourbatov@gmail.com}
\end{center}

\vskip .2 in
\begin{abstract}\noindent
Combining the Hardy-Littlewood $k$-tuple conjecture with a heuristic application of extreme-value statistics,
we propose a family of estimator formulas for predicting maximal gaps between prime $k$-tuples.
Extensive computations show that the estimator $a\log(x/a)-ba$ satisfactorily predicts the maximal gaps below $x$,
in most cases within an error of $\pm2a$, where $a=C_k\log^kx$ is the expected average gap between the same type of $k$-tuples.
Heuristics suggest that maximal gaps between prime $k$-tuples near $x$
are asymptotically equal to $a\log(x/a)$, and thus have the order $O(\log^{k+1}x)$.
The distribution of maximal gaps around the ``trend'' curve $a\log(x/a)$ is close to the Gumbel distribution.
We explore two implications of this model of gaps: record gaps between primes and
Legendre-type conjectures for prime $k$-tuples. 
\end{abstract}

\section{Introduction}
Gaps between consecutive primes have been extensively studied.
The prime number theorem \cite[p.\,10]{hw}
suggests that ``typical'' prime gaps near $p$ have the size about $\log p$. 
On the other hand, maximal prime gaps grow no faster than $O(p^{0.525})$ \cite[p.\,13]{hw}.
Cram\'er \cite{cram} conjectured that gaps between consecutive
primes $p_{n}-p_{n-1}$ are at most about as large as $\log^2p$, that is,
$\lim \sup (p_{n}-p_{n-1})/\log^2p_n = 1$ when $p_n \to \infty$. 
Moreover, Shanks~\cite{shanks} stated that maximal prime gaps $G(p)$ 
satisfy the asymptotic equality $\sqrt{G(p)}\sim\log p$. 
All maximal gaps between primes are now known, 
up to low 19-digit primes (OEIS A005250) \cite{oeis},\,\cite{nicely}. 
This data apparently supports the Cram\'er and Shanks conjectures\footnote{While 
the Shanks conjecture $\sqrt{G(p)}\sim\log p$ is plausible, 
the ``inverted'' Shanks conjecture $p\sim e^{\sqrt{G(p)}}$ is likely false. 
(In general, $X\sim Y \not\Rightarrow e^X\sim e^Y$; for example, $x+\log x\sim x$, 
but $e^{x+\log x}\mspace{-4mu}= xe^x \not\sim e^x$ as $x\to\infty$.)
Wolf \cite[p.\,21]{wolf} proposes an improvement: a gap $G(p)$ is likely to first appear near $p\sim\sqrt{G(p)}e^{\sqrt{G(p)}}$.
}: 
thus far, if we divide by $\log^2p$ the maximal gap ending at $p$,  
the resulting ratio is always less than one~--- but tends to grow 
closer to one, albeit very slowly and irregularly.

Less is known about maximal gaps between prime constellations, or prime $k$-tuples. 
One can conjecture that average gaps between prime $k$-tuples near $p$ 
are $O(\log^{k}p)$ as $p \to \infty$, in agreement with the Hardy-Littlewood $k$-tuple conjecture \cite{hl}.  
Kelly and Pilling \cite{kpII}, 
Fischer~\cite{fischer} 
and Wolf \cite{wolf2} 
report heuristics and computations for gaps between {\em twin primes} ($k=2$).
Kelly and Pilling \cite{kpIII} also provide physically-inspired heuristics for prime triplets ($k=3$);
Fischer \cite{fischer2} conjectures formulas for maximal gaps between $k$-tuples for both $k=2$ and $k=3$.
All of these conjectures and heuristics, as well as extensive computations, suggest that 
maximal gaps between prime $k$-tuples are at most about $\log p$ times the average gap, 
which implies that maximal gaps are $O(\log^{k+1}p)$ as $p \to \infty$. 

In this article we use {\em extreme value statistics} to derive a general formula 
predicting the size of record gaps between $k$-tuples below $p$: 
maximal gaps are approximately 
$a\log(p/a) - ba$, with probable error $O(a)$.
Here $a = C_k\log^kp$ is the {\em expected average gap} near $p$, and 
$C_k$ and $b$ are parameters depending on the type of $k$-tuple. 
This formula approximates maximal gaps 
better and in a wider range than a linear function of $\log^{k+1}p$. 
We will mainly focus on three types of prime $k$-tuples:
\begin{itemize}
\item $k=2$: twin primes (maximal gaps are OEIS A113274);
\item $k=4$: prime quadruplets (maximal gaps are OEIS A113404);
\item $k=6$: prime sextuplets (maximal gaps are OEIS A200503).
\end{itemize}
The observations can be readily applied to other $k$-tuples; 
however, numerical values of constants $C_k$ will change depending on the specific type of $k$-tuple. 
See, e.\,g., the following OEIS sequences for data on maximal gaps between prime $k$-tuples for other $k$: 
\begin{itemize}
\item $k=3$: prime triplets (maximal gaps are A201596 and A201598);
\item $k=5$: prime quintuplets (maximal gaps are A201073 and A201062);
\item $k=7$: prime septuplets (maximal gaps are A201051 and A201251);
\item $k=10$: prime decuplets (maximal gaps are A202281 and A202361). 
\end{itemize}

\section{Definitions, notations, examples}

\noindent
{\em Twin primes} are pairs of consecutive primes that have the form $\{p$, $p+2\}$. 
(This is the densest repeatable pattern of two primes.)
{\em Prime quadruplets} are clusters of four consecutive primes of the form $\{p$, $p+2$, $p+6$, $p+8\}$ 
(densest repeatable pattern of four primes).
{\em Prime sextuplets} are clusters of six consecutive primes of the form $\{p$, $p+4$, $p+6$, $p+10$, $p+12$, $p+16\}$ 
(densest repeatable pattern of six primes).

\noindent
{\em Prime $k$-tuples} are clusters of $k$ consecutive primes that have a repeatable pattern. 
Thus, twin primes are a specific type of prime $k$-tuples, with $k = 2$; 
prime quadruplets are another specific type of prime $k$-tuples, with $k = 4$; 
and prime sextuplets are yet another type of prime $k$-tuples, with $k = 6$. 
(The densest $k$-tuples possible for a given $k$ may also be called 
{\em prime constellations} or {\em prime $k$-tuplets}.)

\noindent
{\em Gaps} between prime $k$-tuples are 
distances between the initial primes in two consecutive $k$-tuples of the same type. 
If the prime at the end of the gap is $p$, we denote the gap $g_k(p)$. 
For example, the gap between the quadruplets $\{11,13,17,19\}$ and $\{101,103,107,109\}$ is $g_4(101) = 90$. 
The gap between the twin primes $\{17,19\}$ and $\{29,31\}$ is $g_2(29) = 12$.
Hereafter $p$ always denotes a prime. 
In the context of gaps between prime $k$-tuples, $p$ will refer to 
the first prime of the $k$-tuple at the end of the gap; we call $p$ the {\em end-of-gap} prime.
(Note that primes {\em preceding} the gap might be {\em orders of magnitude smaller} than the gap size itself;
e.\,g., the gap $g_6(16057)=15960$  starts at $\{97,101,103,107,109,113\}$;
the gap $g_6(1091257)=1047480$ starts at $\{43777,43781,43783,43787,43789,43793\}$.)

\noindent
A {\em maximal gap} is a gap that is strictly greater than all preceding gaps. 
In other words, a maximal gap is the first occurrence of a gap at least this size. 
As an example, consider gaps between prime quadruplets (4-tuples): 
the gap of 90 preceding the quadruplet $\{101,103,107,109\}$ is a maximal gap 
(i.e. the first occurrence of a gap of at least 90), while 
the gap of 90 preceding $\{191,193,197,199\}$ is not a maximal gap 
(not the first occurrence of a gap at least this size). 
A synonym for {\em maximal gap} is {\em record gap}. 
By $G_k(x)$ we will denote the largest gap between $k$-tuples below $x$. 
({\em Note:} Statements like this will always refer to a {\em specific type of $k$-tuples}.)
We readily see that
\begin{eqnarray*}
g_k(p) &\le& G_k(p) \quad\mbox{ wherever $g_k(p)$ is defined, and} \\
g_k(p) &=& G_k(p) \quad\mbox{ if $g_k(p)$ is a maximal gap.}
\end{eqnarray*}
In rare cases, the equality $g_k(p)=G_k(p)$ may also hold for non-maximal gaps $g_k(p)$;
e.\,g., $g_4(191)=G_4(191)=90$ even though the gap $g_4(191)$ is not maximal.

\noindent
The {\it average gap} between $k$-tuples near $x$ 
is denoted $\overline{g_k(x)}$ and defined here as
$$
\overline{g_k(x)}
=
\frac{\displaystyle\sum_{{1\over2}x \,\le\, p-g_k(p) \,<\, p \,\le\, {3\over2}x}\mspace{-48mu}g_k(p)\hphantom{x}}
     {\mspace{-32mu}\displaystyle\sum_{{1\over2}x \,\le\, p-g_k(p) \,<\, p \,\le\, {3\over2}x}\mspace{-48mu}1}   
= 
\frac{\mbox{the sum of all gaps between $k$-tuples with }p\in[{1\over2}x,{3\over2}x]}
     {\mbox{total count of gaps between $k$-tuples with }p\in[{1\over2}x,{3\over2}x]}
\,.
$$
(The value of $\overline{g_k(x)}$ is undefined if there are less than two $k$-tuples with ${1\over2}x \le p \le {3\over2}x$.)

\noindent
The {\it expected average gap} between $k$-tuples near $x$ (for any $x\ge3$) 
is defined formally as $a = a(x) = C_k \log^k x$, 
where the positive coefficient $C_k$ is determined by the type of the $k$-tuples.
(See the {\em Conjectures} section for further details on this.)

\section{Motivation: is a simple linear fit for $G_k(p)$ adequate?}

The first ten or so terms in sequences of record gaps (e.\,g., A113274, A113404, A200503) 
seem to indicate that maximal gaps between $k$-tuples below $p$ grow about as fast as a linear function of $\log^{k+1}p$.
For {\em twin primes} ($k=2$), Rodriguez and Rivera \cite{rr} gave simple linear approximations of record gaps,
while Fischer \cite{fischer2} and Wolf \cite{wolf2} proposed more sophisticated non-linear formulas. 
Why bother with any non-linearity at all?
Let us look at the data. 
Table~1 presents the least-squares zero-intercept trendlines \cite{hays}, \cite{exceltrendline} 
for record gaps between $k$-tuples below $10^{15}$ ($k=2$, 4, 6).

\newpage
\begin{center}TABLE 1 \\Least-squares trendlines for maximal gaps between prime $k$-tuples ($k=2$, 4, 6).\\[0.5em]
\begin{tabular}{cccc}
\hline
    &\multicolumn{3}{c}{\phantom{\fbox{$1^1$}}Trendline equation for maximal gaps between prime $k$-tuples:}    \\
\multicolumn{1}{l}
{End-of-gap prime $p$}&{twin primes}            &{prime quadruplets}      &{prime sextuplets}       \\
                      &{($k=2$; $\xi=\log^3p$)} &{($k=4$; $\xi=\log^5p$)} &{($k=6$; $\xi=\log^7p$)} \\
[0.5ex]\hline
\vphantom{\fbox{$1^1$}}
$1 < p < 10^6$     & $y=0.4576\xi$ &  $y=0.0627\xi$ & $y=0.0016\xi$  \\
$10^6<p< 10^9$     & $y=0.4756\xi$ &  $y=0.1031\xi$ & $y=0.0147\xi$  \\
$10^9<p<10^{12}$   & $y=0.5203\xi$ &  $y=0.1245\xi$ & $y=0.0181\xi$  \\
$10^{12}<p<10^{15}$& $y=0.5628\xi$ &  $y=0.1451\xi$ & $y=0.0249\xi$  \\
\hline
\end{tabular}
\end{center}
Table 1 shows that, for a fixed $k$, record gaps between $k$-tuples farther from zero have a {\em steeper} trendline
(when plotted against $\log^{k+1}p$). This is not a ``one-slope-fits-all'' situation!
There is a good reason to expect that the same tendency holds in general for any $k$:
As we will see in the next sections, there exist {\em curves} that predict the record gap sizes, 
on average, better than any linear function of $\log^{k+1} p$~--- 
and the farther from zero, the steeper are these curves (approaching certain limit values of slope, $C_k$).
Nevertheless, a linear approximation can also be useful; computations and heuristics suggest that
a linear function of $\log^{k+1}p$ can serve as a convenient {\em upper bound} for gaps. 
For example:
{\em Maximal gaps between twin primes are less than} $0.76\log^3 p$.
In what follows, we will combine the Hardy-Littlewood $k$-tuple conjecture with extreme value statistics 
to better predict the sizes of maximal gaps between prime $k$-tuples of any given type,
accounting for their non-linear growth trend.

\section{Conjectures}

In this section we state several conjectures based on plausible heuristics 
and supported by extensive computations.
As far as rigorous proofs are concerned, we do not even know 
whether there are infinitely many $k$-tuples of a given type~--- 
e.\,g., whether there are infinitely many twin primes for $k$ = 2. 
(The famous {\it twin prime conjecture} thus far remains unproven. 
{\it A fortiori} there is no known proof of the more general $k$-tuple conjecture described below.) 

\subsection{The Hardy-Littlewood $k$-tuple conjecture}

The Hardy-Littlewood $k$-tuple conjecture \cite{hl}, \cite[pp.~60--68]{riesel}
predicts the approximate total counts of prime $k$-tuples (with a given admissible\footnote{Any 
pattern of $k$ primes is deemed admissible (repeatable) 
unless it is prohibited by divisibility considerations.
For instance, the pattern of $\{p$, $p+2$, $p+4\}$ is prohibited: 
one of the numbers $p$, $p+2$, $p+4$ must be divisible by 3. 
But $\{p$, $p+2$, $p+6\}$ is not prohibited, hence admissible. 
For a more detailed discussion of admissible patterns, see \cite[pp.\,62--63]{riesel}.
}
pattern):
$$
\mbox{The total number of prime $k$-tuples below $x$} ~\sim~ H_k \int_2^x {dt \over \log^k t}.
$$
The actual counts of $k$-tuples match this prediction with a surprising accuracy \cite[p.\,62]{riesel}.
The coefficients $H_k$ are called the Hardy-Littlewood constants.
Note that, in general, the constants $H_k$ depend on $k$ {\it and} on the specific type of $k$-tuple
(e.\,g., there are three types of prime octuplets, with two different constants).
Hardy and Littlewood not only conjectured the above integral formula
but also provided a recipe for computing the constants $H_k$ as products over subsets of primes.
For example, in special cases with $k=2,4,6$ we have
\begin{eqnarray*}
H_2 &=& 2 \prod_{p\ge3} {{p(p-2)}\over{(p-1)^2}} ~\approx~ 1.32032  \qquad
\mbox{ (for twin primes),} \\
H_4 &=& {27\over2} \prod_{p\ge5} {{p^3(p-4)}\over{(p-1)^4}} ~\approx~ 4.15118 \qquad 
\mbox{ (for prime quadruplets),} \\
%
H_6 &=& {15^5\over2^{13}} \prod_{p\ge7} {{p^5(p-6)}\over{(p-1)^6}} ~\approx~ 17.2986 \qquad 
\mbox{ (for prime sextuplets)}. 
\end{eqnarray*}
These formulas for $H_k$ have slow convergence. 
Riesel \cite{riesel} and Cohen \cite{cohen} describe efficient methods for
computing $H_k$ with a high precision. Forbes \cite{forbes} provides 
the values of $H_k$ for dense $k$-tuples, or {\it $k$-tuplets}, up to $k=24$.
The $k$-tuple conjecture implies that
\begin{itemize}
\item The sequence of maximal gaps between prime $k$-tuples of any given type is infinite.
      (Thus, all OEIS sequences mentioned in {\em Introduction} are infinite.)
\item When $x\to\infty$, the largest gaps below $x$ will grow (asymptotically) at least as fast as average gaps, i.\,e., 
      as fast as $O(\log^{k}x)$ or {\it faster}.
\end{itemize}
But exactly how much faster? Conjectures (D) and (E) below give plausible answers.

\subsection{Conjectured asymptotics for gaps between $k$-tuples}

Let $C_k$ denote the reciprocal to the corresponding Hardy-Littlewood constant: $C_k=H_k^{-1}$.
The following formulas provide rough estimates of the gap $g_k(p)$ ending at a prime $p$: 

\smallskip\noindent
(A) {\it Average} gaps between prime $k$-tuples near $p$ are 
$
\overline{g_k(p)} \sim C_k\log^k p.
$

\smallskip\noindent
(B) {\it Maximal} gaps between prime $k$-tuples are $O(\log^{k+1} p)$: 
$$
g_k(p) < M_k \log^{k+1} p, \quad \mbox{ where } M_k \approx C_k \quad \mbox{ (and possibly } M_k=C_k).
$$
Defining the {\it expected average gap} near $x$ to be
$a = C_k \log^{k} x$~ ($x\ge3$), we further conjecture:

\noindent
(C) {\it Maximal} gaps below $x$ are asymptotically equal to $C_k \log^{k+1} x$:
$$G_k(x) \sim C_k \log^{k+1}x \quad \mbox{ as } x \to \infty, \quad \mbox{with probable error } O(a\log a).
$$

\noindent
(D) {\it Maximal} gaps below $x$ are more accurately described by this asymptotic equality: 
$$
G_k(x) ~\sim~ a\log(x/a) \quad \mbox{ as } x \to \infty, 
\qquad \mbox{ with probable error } O(a).
$$

\noindent
(E) For any given type of $k$-tuple, there exists a real $b$ (e.\,g., $b\approx{2\over k}$) 
such that the difference $G_k(x) - a(\log(x/a)-b)$ changes its sign infinitely often\footnote{ 
Moreover, on finite intervals $x\in[3,X_{\mbox{\tiny max}}]$ the difference $G_k(x) - a(\log(x/a)-b)$ 
changes its sign more often than $G_k(x) - L(\log^{k+1}\negthinspace x)$, 
where $L(\log^{k+1}\negthinspace x)$ is any linear function of $\log^{k+1}\negthinspace x$ 
and $X_{\mbox{\tiny max}}$ is large enough.
} 
as $x\to\infty$.

\smallskip
A key ingredient in these conjectures is provided by the constants $C_k = {H_k}^{-1}$: 
$$
C_2 = H_2^{-1} \approx 0.75739,  \qquad
C_4 = H_4^{-1} \approx 0.240895, \qquad
C_6 = H_6^{-1} \approx 0.057808.
$$
Another key ingredient is a statistical formula: for certain kinds of random events occurring
%
%
at mean intervals $a$, the record interval between events observed in time $T$ is likely\footnote{
In particular, if intervals between rare random events have the exponential distribution, with mean interval $a$ sec
and CDF $1-e^{-t/a}$, then the {\bf\em most probable record interval} 
observed within $T$ sec is about $a\log(T/a)$ sec (provided that $a\ll T$).
After many observations ending at times $a \ll T_1 \ll T_2 \ll T_3\ldots$ 
almost surely for some $T_i$ we will observe record intervals {\bf\em exceeding} $a\log(T_i/a)$.
However, for other values of $T_i$ we will also observe record intervals {\bf\em below} $a\log(T_i/a)$.
It is this formula for the most probable extreme, with the aid of the estimate SD\,$=O(a)$ 
for the standard deviation of extremes, that allows us to heuristically predict the bounds, errors, asymptotics, 
and sign changes in conjectures (B), (C), (D), (E).
}
near $a\log(T/a)$.
In {\em Appendix} we derive this formula for $a=\mbox{const}$. 
Here, we heuristically apply this formula for a slowly changing $a$ (i.\,e., $a = C_k \log^{k} x$).
For now, we can informally summarize the behavior of maximal gaps between $k$-tuples near $p$ as follows:
Maximal gaps are {\bf{\em at most}} about $\log p$ times the average gap.

\subsection{Estimators for maximal gaps between $k$-tuples}

Prime $k$-tuples are rare and seemingly ``random''. 
Life offers many examples of unusually large intervals between rare random events, 
such as the longest runs of dice rolls without getting a twelve; 
maximal intervals between clicks of a Geiger counter measuring very low radioactivity, etc. 
Reasoning as in {\em Appendix}, one can statistically estimate the mathematical expectation of maximal intervals 
between rare random events by expressing them in terms of the average intervals: 
$$
\mbox{Expected maximal intervals} ~=~ a \log (T/a) + O(a), \eqno{(*)}
$$
where $a$ is the average interval between the rare events, and 
$T$ is the total observation time or length ($1 \ll a \ll T$). 

To account for the observed non-linear growth of record gaps  between prime $k$-tuples (Table 1), 
we will simulate gap sizes using estimator formulas very similar to the above $(*)$.
We define a family of estimators for the maximal gap that ends at $p$:
\begin{eqnarray}
E_1(G_k(p)) &=& \max(a,\, a\log(p/a) - ba),      \qquad    \mbox{ probable error: $O(a);$ }\\
E_2(G_k(p)) &=& \max(a,\, a\log(p/a)),    \qquad\qquad\,   \mbox{ probable error: $O(a);$ }\\
E_3(G_k(p)) &=& a\log p ~=~ C_k \log^{k+1} p, \qquad\quad  \mbox{ probable error: $O(a\log a).$ }
\end{eqnarray}
Here, the role of the statistically average interval $a$ is played by the expected average gap between $k$-tuples: 
as before, we set $a = C_k \log^k p$. 
The role of the total observation time $T$ is played by $p$ (we are ``observing'' gaps that occur from 0 to $p$).
We also empirically choose $b={2\over k}$.
(The latter choice is not set in stone; by varying the parameter $b$ in $E_1$ one can get 
an infinite family of useful estimators
with similar asymptotics. In Section 6 we will see that $b\approx3$ appears
quite suitable for modeling {\em prime gaps}, in which case $k=1$, $a=\log p$, and $C_1=1$.) 
It is easy to see that, for any fixed $k\ge1$ 
and any fixed $b\ge0$, we have
$$
a \le E_1 \le E_2 < E_3 \mbox{ ~for all } p\ge3, \quad \mbox{ but at the same time } \quad 
a \ll E_1 \sim E_2 \sim E_3 \mbox{ ~as } p \to \infty.
$$
Indeed, when $p \to \infty$ we have
$$
E_1(G_k(p)) = a \log(p/a)-ab = a\log p-a(\log a + b) = a\log p- o(\log^{k+1}p) \sim E_3(G_k(p)).
$$
{\em Note:} We use the max function in the estimators to guarantee that $E_i(G_k(p)) \ge a$.
This precaution is needed because, if $p$ is not large enough, $\log(p/a)$ might be negative or too small.
We want our estimators to give positive predictions no less than $a$ even in such cases.

\medskip\noindent
The above conjectures (C), (D), (E) tell us that $E_1$ and $E_2$ are better estimators than $E_3$:
the probable error of $E_3$ is greater than that of $E_1$ or $E_2$. 
In Section 5.1 we will compare the predictions obtained with these estimators 
to the actual sizes of maximal gaps. 

\subsection{Why extreme value statistics?}

In number theory, probabilistic models such as Cram\'er's model \cite{cram} face serious difficulties.
One such difficulty will be noted in Section 6.
Pintz \cite{cvc} points out additional problems with Cram\'er's model.
Number-theoretic objects (such as primes or prime $k$-tuples) are too peculiar; 
they are clearly {\em not} independent and cannot be flawlessly simulated by 
independent and identically distributed (i.i.d.)~random variables 
or ``events'' or ``coin tosses'' that we usually deal with in probabilistic models.
Why then should one build heuristics for prime $k$-tuples based on extreme value statistics?

An obvious reason is that we are studying {\em extreme} gaps, so it would be unwise 
to outright dismiss the existing extreme value theory without giving it a try.
%
%
When our goal is just to guess the right formula, rigor is not the highest priority; 
it is perhaps more important to accumulate as much evidence as possible,
look for counterexamples, and make reasonable simplifications. 
The above formula for the expected maximal interval $(*)$ appears to be at the right level of simplification
and fits the actual record gaps fairly well even without the $O(a)$ term (as we will see in Section 5.1).
To fine-tune formula $(*)$ for record gaps between prime $k$-tuples, we simply have to find a suitable $O(a)$ term. 
The latter can be done using number-theoretic insights and/or numerical evidence.

Extreme value theory also offers additional benefits.
Not only does it tell us the mathematical expectation of extremes in random sequences~--- 
it also predicts {\it distributions} of extremes.
While in general there are infinitely many probability distribution laws,
there exist {\em only three types of limiting extreme value distributions}
applicable to sequences of i.i.d.~random variables: 
the Gumbel, Fr\'echet, and Weibull distributions \cite{bgt}.
%
%
When no limiting extreme value distribution exists, a known type of 
extreme value distribution may still be a good approximation 
\cite{gsw},\,\cite{schilling}.
A large body of knowledge has been accumulated that extends the same types of extreme value distributions 
from i.i.d.~random variables to certain kinds of dependent variables, for example, 
$m$-dependent random variables \cite{watson},
exchangeable variables \cite{berman},\,\cite[pp.\,163--191]{galambos}, 
and other situations \cite{bgt},\,\cite{galambos}.
Although no theorem currently extends the known types of extreme value distributions 
to record gaps between primes or prime $k$-tuples, 
we might have an aesthetic expectation 
that ``the usual suspects'' would show up here, too.
It turns out that one common type of extreme value distribution~--- the Gumbel distribution~--- does show up! 
(See Section 5.2, {\em The distribution of maximal gaps}.)

\pagebreak\newpage

\section{Numerical results}

Using a fast deterministic algorithm based on strong pseudoprime tests \cite[pp.~91--92]{riesel}, 
the~author computed all maximal gaps between prime $k$-tuplets up to $10^{15}$ for $k=4$,~6. 
Fischer (2008) \cite{fischer} reported a similar computation for $k=2$. Below we analyze this data.


\subsection{The growth of maximal gaps}

Figure 1 shows record gaps between twin primes (A113274) for $p<10^{15}$;
the curves are predictions obtained with estimators $E_1$, $E_2$, $E_3$ defined above.
Figure 2 shows similar data for prime quadruplets (A113404),
and Figure 3 for prime sextuplets (A200503).
Tables 2--4 give the relevant numerical data; see also OEIS sequences mentioned in {\em Introduction}.

\begin{figure}[hbtp] 
  \centering
  \includegraphics[bb=21 38 765 513,width=5.67in,height=5.67in,keepaspectratio]{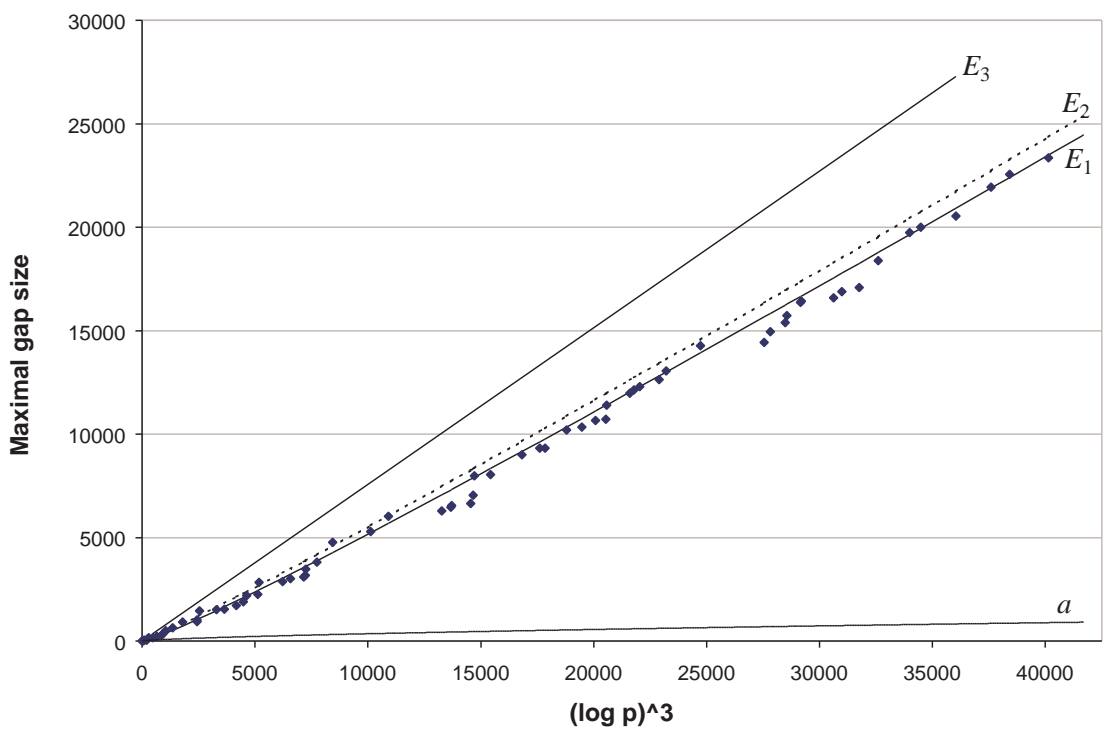}
  \caption{Maximal gaps between twin primes $\{p, p+2\}$ 
           (A113274).
           Plotted (bottom to top): expected average gap $a=0.75739\log^2 p$, 
           estimators $E_1=a\log(p/a)-ba$, $E_2=a\log(p/a)$, $E_3=a\log p=0.75739\log^3 p$,
           where $p$ is the end-of-gap prime; $b=1$.
  }
  \label{fig:2tuples-maximal-gaps}
\end{figure}

Here are some observations suggested by these numerical results. 
(As before, $a$ denotes the expected average gap, $a=C_k \log^k p$, and $b={2\over k}$ unless stated otherwise.)

\begin{enumerate}
\item Estimators $E_1$ and $E_2$ overestimate some of the actual record gaps, but underestimate others. 
      For $k\le6$, the data shows that $E_1$ is closer to a median-unbiased
      estimator.\footnote{A median-unbiased estimator $E_{\mbox{\tiny med}}(x)$ has as many observed values above it as below it.}
      (We can make it even closer by tweaking the $b$ value; 
      e.\,g., setting $b\approx1.2597$ for twin primes, or $b\approx0.7497$ for prime quadruplets, 
      would turn $E_1$ into a median-unbiased estimator for maximal gaps below $10^{15}$.)
\item About 90\% of the observed gaps are within $\pm2a$ of the~$E_1$ curve.
      Over 50\% of the observed gaps are within $\pm a$ of $E_1$. 
This level of accuracy appears to be in line with heuristics based on statistical models 
(where the relevant extreme-value distributions have 
the standard deviation $\pi a/\sqrt{6}\approx1.28a$; see {\em Appendix}).
\item Consider median-unbiased estimators $E_{\mbox{\tiny med}}(G_k(p))=a(\log(p/a)-b_{\mbox{\tiny med}})$ for $p<10^{15}$.
      Computations show that the value of $b_{\mbox{\tiny med}}$ tends to decrease when $k$ increases; also,
      our empirical value $b={2\over k}$ in the $E_1$ estimator 
      is a little closer to zero than the median-unbiased value $b_{\mbox{\tiny med}}$.
      (For a simple way to refine $b$, see remark at the end of sect.\,5.2.)
\item For relatively small values of $p$ that we deal with, the estimator $E_3$ may seem useless (too far above the realistic values). 
      However, all three estimators are asymptotically equivalent,
      $E_1 \sim E_2 \sim E_3$ when $p\to\infty$.

\item The estimator $E_3=C_k \log^{k+1}p$ overestimates all known record gaps. 
      In most cases, the error of $E_3$ is close to $a \log a$, exactly as expected from extreme-value statistics.
      Thus $E_3$ may be a good candidate for an {\em upper bound} for all record gaps;
      so in statement (B) of section 4.2 we may have $M_k=C_k=H_k^{-1}$, and
$$
G_k(p)<C_k \log^{k+1}p \qquad\mbox{(an analog of Cram\'er's conjecture)}.
$$
It would be interesting to see any counterexamples, i.\,e., gaps exceeding $C_k \log^{k+1}p$.
\begin{figure}[hb] 
  \centering
  \includegraphics[bb=22 15 780 495,width=5.67in,height=5.67in,keepaspectratio]{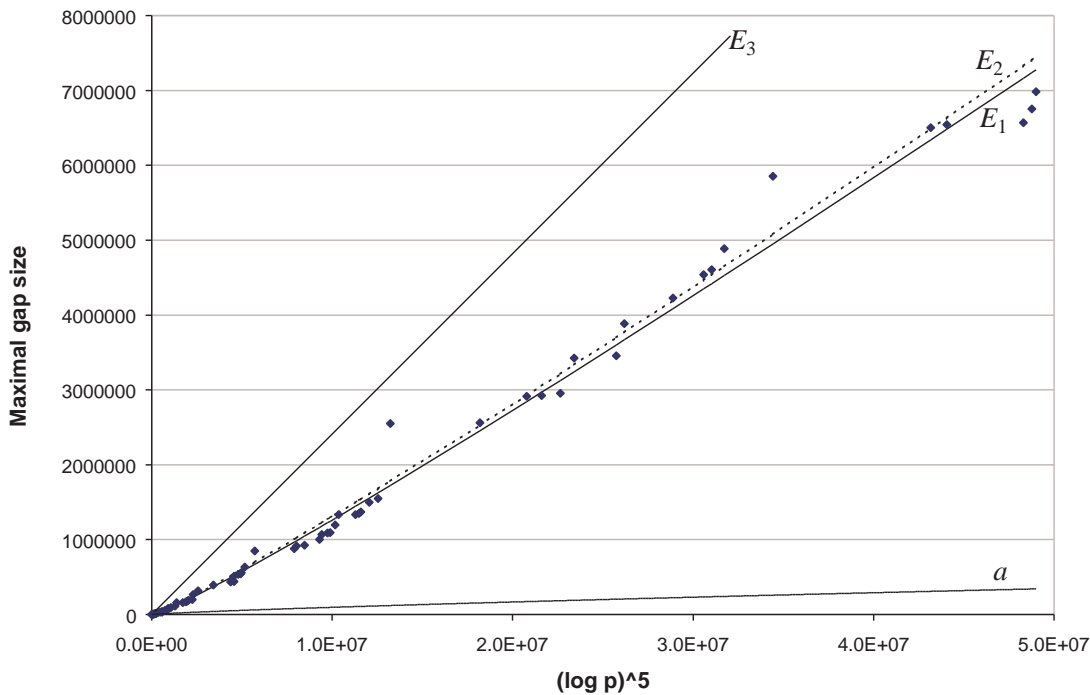}
  \caption{Maximal gaps between prime quadruplets
           $\{p, p+2, p+6, p+8\}$ 
           (A113404).
           Plotted (bottom to top): expected average gap $a=0.24089\log^4 p$, 
           estimators $E_1=a\log(p/a)-ba$, $E_2=a\log(p/a)$, $E_3=a\log p=0.24089\log^5 p$,
           where $p$ is the end-of-gap prime; $b=1/2$.
  }
  \label{fig:4tuples-maximal-gaps-1204a}
\end{figure}
\end{enumerate}
 
\smallskip\noindent
{\bf Absolute error.} The absolute error $|E_i-G_k(p)|$ tends to grow (but not monotonically) as $p\to\infty$
for all three estimators $E_1$, $E_2$, $E_3$.
Heuristically, we expect the absolute error to be unbounded and, on average, continue to grow 
for all three estimators. Probable absolute errors are $O(a)$ for $E_1$ and $E_2$, and $O(a\log a)$ for $E_3$. 

\begin{figure}[htb] 
  \centering
  \includegraphics[bb=13 11 783 500,width=5.67in,height=5.67in,keepaspectratio]{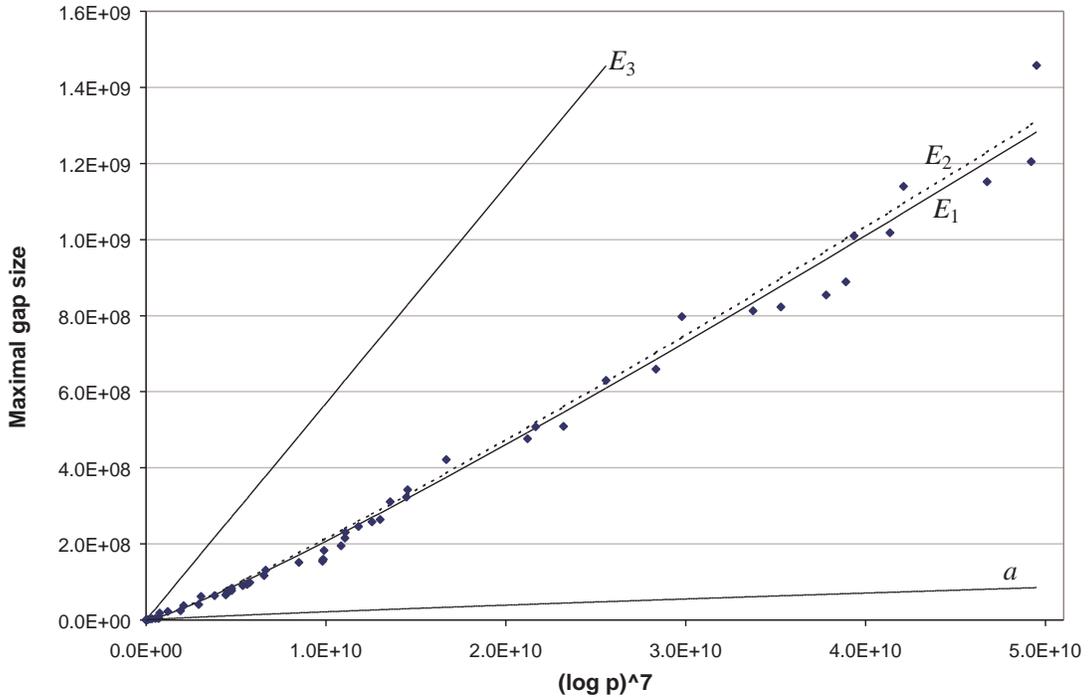}
  \caption{Maximal gaps between prime sextuplets $\{p$, $p+4$, $p+6$, $p+10$, $p+12$, $p+16\}$ 
           (A200503).
           Plotted (bottom to top): expected average gap $a=0.057808\log^6 p$, 
           estimators $E_1=a\log(p/a)-ba$, $E_2=a\log(p/a)$, $E_3=a\log p=0.057808\log^7 p$,
           where $p$ is the end-of-gap prime; $b=1/3$. 
  }
  \label{fig:6tuples-maximal-gaps}
\end{figure}

\smallskip\noindent
{\bf Relative error.} The relative error $|\varepsilon_i|=|E_i-G_k(p)|/G_k(p)$ tends to decrease (but not monotonically) 
for all three estimators as $p\to\infty$.
It may not be obvious from Figures~\ref{fig:2tuples-maximal-gaps}--\ref{fig:6tuples-maximal-gaps}, 
but we must have $|\varepsilon_i|\to0$ either for {\it all three} estimators or for {\it none} of them.
({\em Note}: the limit of $|\varepsilon_i|$ as $p\to\infty$ might not exist at all; 
that would invalidate most of our conjectures.)

\smallskip\noindent
{\bf Error in average-gap units {\em a}.} 
The error $(E_3-G_k(p))/a$, i.\,e., the $E_3$ error expressed as a number of expected average gaps, 
grows about as fast as $\log a$ (but not monotonically). 
Judging from limited numerical data, the corresponding error $(E_i-G_k(p))/a$ seems bounded as $p\to\infty$ 
if we use estimators $E_1$ or $E_2$. 
Heuristically, for $E_1$ and $E_2$ this error should remain bounded for the majority (but not all) of the record gaps.

\smallskip\noindent
Overall, the 
prediction that record gaps are about 
$a\log(p/a)+O(a)$ appears correct for the vast majority of actual gaps, as far as we have checked ($p<10^{15}$). 
Note that the ``optimal'' $O(a)$ term ($-ba$ in the $E_1$ estimator) is negative, at least for $k\le6$.
For larger values of $k$, the parameter $b$ gets closer to zero. 
Empirically, for $k$-tuples with $k\ge6$, 
the $E_1$ estimator will likely produce good results even with $b\approx0$. 
Therefore, for large $k$ we might want to simplify the model and use $b=0$, i.\,e., use the estimator $E_2=\max(a,a\log(p/a))$, 
the dotted curve in the above figure. 
However, maximal gap estimators with certain special properties 
(e.g., median-unbiased estimators) will still require nonzero values of $b$.

\phantom{\vspace{5mm}}

\begin{center}TABLE 2 \\Maximal gaps between twin primes $\{p$, $p+2\}$ below $10^{15}$ \\[0.5em]
\begin{tabular}{rrr|rrr}  
\hline
End-of-gap prime  & Gap $g_2$ & $g^*_2$~\vphantom{\fbox{$1^1$}} & 
End-of-gap prime  & Gap $g_2$ & $g^*_2$~\vphantom{\fbox{$1^1$}} \\
[0.5ex]\hline
\vphantom{\fbox{$1^1$}}
5&2&0.084&24857585369&6552&$-2.765$ \\
11&6&0.451&40253424707&6648&$-3.585$\\
29&12&0.180&42441722537&7050&$-2.811$\\
59&18&$-0.115$&43725670601&7980&$-0.830$\\
101&30&0.025&65095739789&8040&$-1.625$\\
347&36&$-1.205$&134037430661&8994&$-1.323$\\
419&72&$-0.112$&198311695061&9312&$-1.604$\\
809&150&1.247&223093069049&9318&$-1.865$\\
2549&168&$-0.396$&353503447439&10200&$-1.262$\\
6089&210&$-1.011$&484797813587&10338&$-1.747$\\
13679&282&$-1.189$&638432386859&10668&$-1.792$\\
18911&372&$-0.486$&784468525931&10710&$-2.195$\\
24917&498&0.645&794623910657&11388&$-1.032$\\
62927&630&0.290&1246446383771&11982&$-1.081$\\
188831&924&0.834&1344856603427&12138&$-0.998$\\
688451&930&$-1.728$&1496875698749&12288&$-1.002$\\
689459&1008&$-1.161$&2156652280241&12630&$-1.309$\\
851801&1452&1.577&2435613767159&13050&$-0.916$\\
2870471&1512&$-0.721$&4491437017589&14262&$-0.481$\\
4871441&1530&$-1.689$&13104143183687&14436&$-2.773$\\
9925709&1722&$-2.070$&14437327553219&14952&$-2.255$\\
14658419&1902&$-1.948$&18306891202907&15396&$-2.181$\\
17384669&2190&$-0.918$&18853633240931&15720&$-1.793$\\
30754487&2256&$-1.805$&     23275487681261 &16362&$-1.398$\\
32825201&2832&0.601  &     23634280603289 &16422&$-1.351$\\
96896909&2868&$-1.646$&     38533601847617 &16590&$-2.291$\\
136286441&3012&$-1.812$&   43697538408287 &16896&$-2.178$\\
234970031&3102&$-2.611$&   56484333994001 &17082&$-2.539$\\
248644217&3180&$-2.451$&   74668675834661 &18384&$-1.507$\\
255953429&3480&$-1.454$&  116741875918727&19746&$-0.864$\\
390821531&3804&$-1.260$&  136391104748621&19992&$-0.940$\\
698547257&4770&0.571&    221346439686641&20532&$-1.467$\\
2466646361&5292&$-0.816$& 353971046725277&21930&$-0.955$\\
4289391551&6030&$-0.075$& 450811253565767&22548&$-0.834$\\
19181742551&6282&$-2.831$&742914612279527&23358&$-1.149$\\
24215103971&6474&$-2.888$& & & \\
\hline
\end{tabular}
\end{center}
\normalsize

{\vskip 2mm}
\noindent
%

\newpage

\begin{center}TABLE 3 \\Maximal gaps between prime quadruplets $\{p$, $p+2$, $p+6$, $p+8\}$ below $10^{15}$ \\[0.5em]
\begin{tabular}{rrr|rrr}  
\hline
End-of-gap prime  & Gap $g_4$ & $g^*_4$~\vphantom{\fbox{$1^1$}} & 
End-of-gap prime  & Gap $g_4$ & $g^*_4$~\vphantom{\fbox{$1^1$}} \\
[0.5ex]\hline
\vphantom{\fbox{$1^1$}}
11&6&0.430&3043668371&557340&$-0.750$\\
101&90&0.902&3593956781&635130&0.188\\
821&630&0.770&5676488561&846060&2.366\\
1481&660&0.192&25347516191&880530&$-1.576$\\
3251&1170&$-0.014$&27330084401&914250&$-1.358$\\
5651&2190&0.194&35644302761&922860&$-1.966$\\
9431&3780&0.518&56391153821&1004190&$-2.244$\\
31721&6420&$-0.125$&60369756611&1070490&$-1.697$\\
43781&8940&0.211&71336662541&1087410&$-1.967$\\
97841&9030&$-0.998$&76429066451&1093350&$-2.089$\\
135461&13260&$-0.539$&87996794651&1198260&$-1.383$\\
187631&16470&$-0.434$&96618108401&1336440&$-0.242$\\
326141&24150&$-0.094$&151024686971&1336470&$-1.535$\\
768191&28800&$-1.004$&164551739111&1348440&$-1.663$\\
1440581&29610&$-1.957$&171579255431&1370250&$-1.577$\\
1508621&39990&$-0.977$&211001269931&1499940&$-0.986$\\
3047411&56580&$-0.812$&260523870281&1550640&$-1.150$\\
3798071&56910&$-1.217$&342614346161&2550750&6.412\\
5146481&71610&$-0.714$&1970590230311&2561790&0.197\\
5610461&83460&$-0.049$&4231591019861&2915940&$-0.076$\\
9020981&94530&$-0.379$&5314238192771&2924040&$-0.748$\\
17301041&114450&$-0.678$&7002443749661&2955660&$-1.421$\\
22030271&157830&0.996&8547354997451&3422490&0.447\\
47774891&159060&$-0.861$&15114111476741&3456720&$-1.200$\\
66885851&171180&$-1.135$&16837637203481&3884670&0.533\\
76562021&177360&$-1.202$&30709979806601&4228350&0.134\\
87797861&190500&$-1.023$&43785656428091&4537920&0.307\\
122231111&197910&$-1.515$&47998985015621&4603410&0.278\\
132842111&268050&0.677&55341133421591&4884900&0.972\\
204651611&315840&1.022&92944033332041&5851320&2.995\\
628641701&395520&0.099&412724567171921&6499710&0.021\\
1749878981&435240&$-1.669$&473020896922661&6544740&$-0.293$\\
2115824561&440910&$-2.020$&885441684455891&6568590&$-2.253$\\
2128859981&513570&$-0.617$&947465694532961&6750330&$-1.932$\\
2625702551&536010&$-0.749$&979876644811451&6983730&$-1.356$\\
2933475731&539310&$-0.982$& & & \\			
\hline
\end{tabular}
\end{center}
\normalsize

{\vskip 2mm}
\noindent
{\em Notes:} Computing Tables 3 and 4 took the author two weeks on a quad-core 2.5 GHz CPU.
Table 2 reflects Fischer's extensive computation \cite{fischer}. 
For earlier computations of maximal gaps by Boncompagni, Rodriguez, and Rivera, see OEIS A113274, A113404 \cite{oeis}, \cite{rr}. 
\normalsize

\newpage

\begin{center}TABLE 4 \\Maximal gaps between prime 6-tuples $\{p$, $p+4$, $p+6$, $p+10$, $p+12$, $p+16\}$ below $10^{15}$ \\[0.5em]
\small
\begin{tabular}{rrr|rrr}  
\hline
End-of-gap prime  & Gap $g_6$ & $g^*_6$~\vphantom{\fbox{$1^1$}} & 
End-of-gap prime  & Gap $g_6$ & $g^*_6$~\vphantom{\fbox{$1^1$}} \\
[0.5ex]\hline
\vphantom{\fbox{$1^1$}}
97&90&1.856&422248594837&159663630&$-2.389$\\
16057&15960&1.414&427372467157&182378280&$-1.353$\\
43777&24360&0.949&610613084437&194658240&$-1.751$\\
1091257&1047480&1.570&660044815597&215261760&$-1.079$\\
6005887&2605680&1.176&661094353807&230683530&$-0.427$\\
14520547&2856000&$-0.048$&853878823867&245336910&$-0.573$\\
40660717&3605070&$-1.035$&1089869218717&258121710&$-0.786$\\
87423097&4438560&$-1.646$&1248116512537&263737740&$-0.968$\\
94752727&5268900&$-1.380$&1475318162947&311017980&0.246\\
112710877&17958150&3.778&1914657823357&322552230&$-0.170$\\
403629757&21955290&1.526&1954234803877&342447210&0.436\\
1593658597&23910600&$-1.149$&3428617938787&421877610&1.085\\
2057241997&37284660&0.730&9368397372277&475997340&$-0.740$\\
5933145847&40198200&$-1.318$&10255307592697&507945690&$-0.256$\\
6860027887&62438460&1.224&13787085608827&509301870&$-1.159$\\
14112464617&64094520&$-0.506$&21017344353277&629540730&0.084\\
23504713147&66134250&$-1.523$&33158448531067&659616720&$-0.819$\\
24720149677&70590030&$-1.228$&41349374379487&797025180&0.985\\
29715350377&77649390&$-1.038$&72703333384387&813158010&$-0.682$\\
29952516817&83360970&$-0.556$&89166606828697&823158840&$-1.190$\\
45645253597&90070470&$-1.064$&122421855415957&854569590&$-1.723$\\
53086708387&93143820&$-1.202$&139865086163977&888931050&$-1.642$\\
58528934197&98228130&$-1.063$&147694869231727&1010092650&$-0.077$\\
93495691687&117164040&$-0.935$&186010652137897&1018139850&$-0.755$\\
97367556817&131312160&$-0.108$&202608270995227&1139590200&0.603\\
240216429907&151078830&$-1.388$&332397997564807&1152229260&$-0.967$\\
414129003637&154904820&$-2.566$&424682861904937&1204960680&$-1.155$\\
419481585697&158654580&$-2.420$&437805730243237&1457965740&1.725\\
\hline
\end{tabular}
\end{center}
\normalsize

\phantom{\vspace{1cm}}

\subsection{The distribution of maximal gaps}

We have just seen that maximal gaps between prime $k$-tuples below $p$ grow about as fast as $a\log(p/a)$.
Thus, the curve $a\log(p/a)$ (the dotted curve in Figures~1--3) may be regarded as a ``trend.''
Now we are going to take a closer look at the distribution of maximal gaps in the neighborhood of this ``trend'' curve. 
In our analysis, we will also include the case $k=1$, record gaps between primes (A005250). 
For each $k=1$, $2$, $4$, $6$, we will make a histogram of shifted and scaled ({\em standardized}) record gaps:
subtract the ``trend'' $a\log(p/a)$ from actual gaps, and then 
divide the result by the ``natural unit'' $a$, the expected average gap.  
This way, all record gaps $g_k(p)$ are mapped to standardized values $g_k^*$ (shown in Tables 2--4):
$$
g_k(p) ~~\to~~ g_k^* = \frac{g_k(p)-a\log(p/a)}{a} 
, \quad \mbox{ where~ } a = C_k \log^k p.
$$
Record gaps that exceed $a\log(p/a)$ are mapped to standardized values $g_k^* > 0$, 
while those below $a\log(p/a)$ are mapped to $g_k^* < 0$.
Note that the majority of known record gaps are below the dotted curve in Figures 1--3; 
accordingly, most of the standardized values $g_k^*$ are negative.
It is also immediately apparent that the histograms and fitting distributions are skewed: 
the right tail is longer and heavier.
This skewness is a well-known characteristic of {\em extreme value distributions}~---
and it comes as no surprise that a good fit obtained with the help of distribution-fitting software \cite{easyfit} 
is the {\em Gumbel distribution}, a common type of extreme value distribution (see {\em Appendix}).
\begin{figure}[tbh] 
\centering
\includegraphics[bb=11 25 778 320,width=6in,height=4.8in,keepaspectratio]{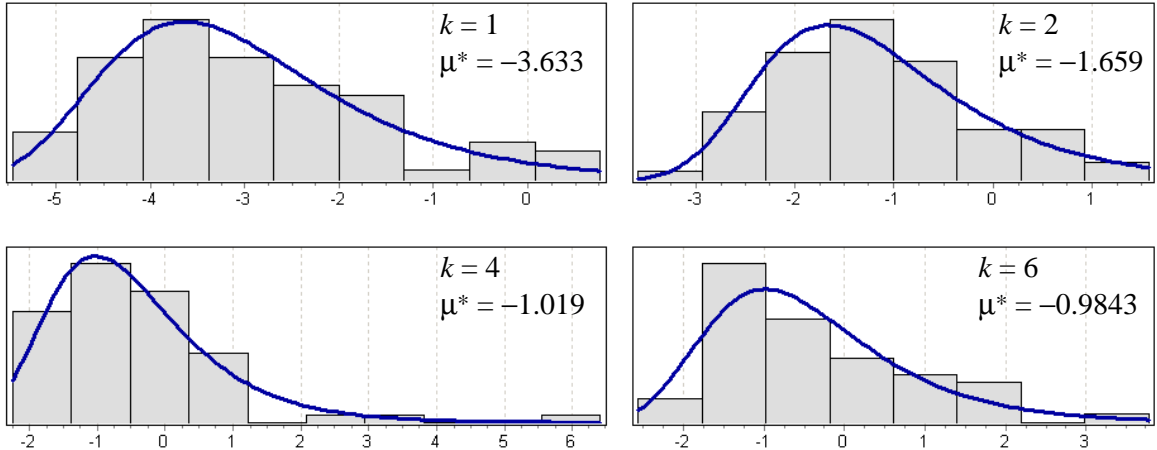}
\caption{The distribution of standardized maximal gaps $g_k^*$: 
histograms and the fitting Gumbel distribution PDFs.
For $k=1$ (primes), the histogram shows record gaps below $4\times10^{18}$.
For $k=2,4,6$ ($k$-tuples), the histograms show record gaps below $10^{15}$. 
}
\label{fig:fig4histograms}
\end{figure}

Here is why we can say that the Gumbel distribution is indeed a good fit:

\smallskip\noindent 
(1) Based on goodness-of-fit statistics (the Anderson-Darling test as well as the Kolmogorov-Smirnov test), 
one {\bf\em cannot reject} the hypothesis that the standardized values 
$g_k^*$ might be values of independent identically distributed random variables with the Gumbel distribution.

\smallskip\noindent
(2) Although a few other distributions could not be rejected either,
the Anderson-Darling and Kolmogorov-Smirnov goodness-of-fit statistics for the Gumbel distribution are better
than the respective statistics for any other two-parameter distribution we tried 
(including~normal, uniform, logistic, Laplace, Cauchy, power-law, etc.),
and better than for several three-parameter distributions 
(e.\,g.,~triangular, error, Beta-PERT, and others).

\smallskip
An equally good or even marginally better fit is the three-parameter
{\em generalized extreme value} (GEV) distribution, which in fact includes the Gumbel distribution as a special case.
The {\em shape parameter} in the fitted GEV distribution turns out very close to zero; 
note that a GEV distribution with a zero shape parameter is precisely the Gumbel distribution.
The {\em scale} parameter of the fitted Gumbel distribution is close to one.
The {\em mode} $\mu^{*}$ of the fitted distribution is negative.
Figure~\ref{fig:fig4histograms} gives the approximate value of $\mu^{*}$ for $k=1,2,4,6$;
$\mu^{*}$ is the coordinate of the maximum of the distribution PDF (probability density function).

{\it Note}: Now that we have a more precise value of the mode $\mu^*$, 
we can refine the parameter $b$ in the $E_1$ estimator: 
use $-b=\mu^* + \gamma$, 
which estimates the {\em mean} of the fitted Gumbel distribution in Fig.\,\ref{fig:fig4histograms}.
Here $\gamma = 0.5772...$ is the Euler-Mascheroni constant.

\newpage
\section{On maximal gaps between primes}

\begin{figure}[h] 
  \centering
  \includegraphics[bb=14 19 783 486,width=5.67in,height=5.67in,keepaspectratio]{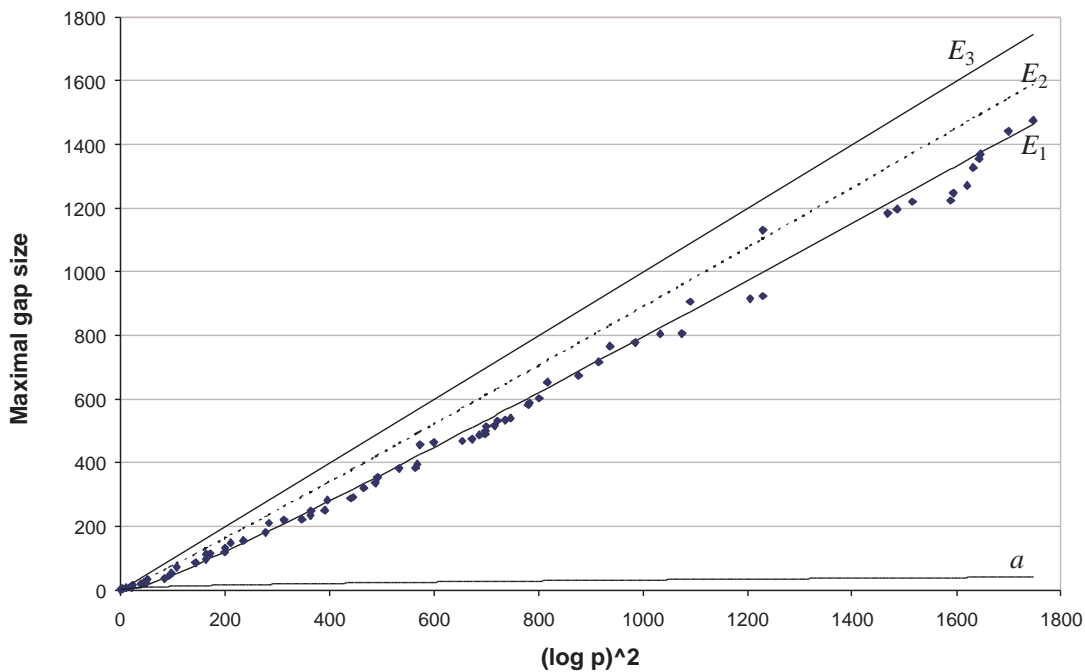}
  \caption{Maximal gaps between consecutive primes (A005250).
           Plotted (bottom to top): expected average gap $a=\log p$, 
           estimators $E_1=a\log(p/a)-ba$, $E_2=a\log(p/a)$ (dotted), $E_3=a\log p=\log^2 p$,
           where $p$ is the end-of-gap prime; $b=3$. 
  }
  \label{fig:1tuples-maximal-gaps}
\end{figure}

Let us now apply our model of gaps to {\em maximal gaps between primes} (A005250) \cite{oeis},\,\cite{nicely}: 
\begin{quotation}
Maximal prime gaps are about $a\log(p/a)-ba$, with $a=\log p$ and $b\approx3$. 
\end{quotation}
If all record gaps behave like those in Figure~\ref{fig:1tuples-maximal-gaps}
(showing the 75 known record gaps between primes $p<4\times10^{18}$), 
this would confirm the Cram\'er and Shanks conjectures: 
maximal prime gaps are smaller than $\log^2 p$~--- but smaller only by $O(a\log a)$. 
We also easily see that the Cram\'er and Shanks conjectures are compatible with our estimate of record gaps.
Indeed, for $a = \log p$ and {\em any fixed} $b\ge0$, we have $\log^2 p\sim a(\log(p/a)-b)<\log^2 p$ as $p\to\infty$. 

\noindent
{\em Notes}:
Maier's theorem (1985) \cite{maier} states that there are (relatively short) intervals 
where typical gaps between primes are greater than the average ($\log p$) expected from the prime number theorem. 
Based in part on Maier's theorem, Granville \cite{granville} adjusted the Cram\'er conjecture and proposed that, 
as $p\to\infty$, $\lim \sup(G(p)/\log^2p) \ge 2e^{-\gamma} = 1.1229\ldots$ 
This would mean that an infinite subsequence of maximal gaps must lie above the Cram\'er-Shanks upper limit 
$\log^2p$, i.\,e., above the $E_3$ line in Figure~\ref{fig:1tuples-maximal-gaps} --- 
and this hypothetical subsequence (or an infinite subset thereof) 
must approach a line whose slope is about 1.1229 times steeper! 
However, for now, there are no known maximal prime gaps above $\log^2p$. 
Interestingly, Maier himself did not voice serious concerns that the Cram\'er or Shanks conjecture might be 
in danger because of his theorem; thus, Maier and Pomerance \cite{mp} simply remarked in 1990:
\begin{quotation}\noindent{\small
Cram\'er conjectured that $\lim\sup G(x)/\log^2\negthinspace x = 1$, while Shanks made the stronger conjecture that 
$G(x)\sim\log^2x$, but we are still a long way from proving these statements.
}
\end{quotation}
%

\section{Corollaries: Legendre-type conjectures}

Assuming the conjectures of Section 4, 
one can state (and verify with the aid of a computer) 
a number of interesting corollaries. 
The following conjectures generalize Legendre's conjecture about primes between squares.
\begin{itemize}
\item For each integer $n > 0$, there is always a prime between $n^2$ and $(n+1)^2$. ({\em Legendre})  
\item For each integer $n > 122$, there are twin primes between $n^2$ and $(n+1)^2$. (A091592) 
\item For each integer $n > 3113$, there is a prime triplet between $n^2$ and $(n+1)^2$. 
\item For each integer $n > 719377$, there is a prime quadruplet between $n^2$ and $(n+1)^2$. 
\item For each integer $n > 15467683$, there is a prime quintuplet between $n^2$ and $(n+1)^2$. 
\item There exists a sequence $\{s_k\}$ such that, for each integer $n > s_k$, 
      there is a prime $k$-tuplet between $n^2$ and $(n+1)^2$. 
      (This $\{s_k\}$ is OEIS A192870: $0, 122, 3113, 719377, \ldots$)
\end{itemize}
Another family of Legendre-type conjectures for prime $k$-tuplets 
can be obtained by replacing squares with cubes, 4th, 5th, and higher powers of $n$:
\begin{itemize}
\item For each integer $n > 0$, there are twin primes between $n^3$ and $(n+1)^3$.
\item For each integer $n > 0$, there is a prime triplet between $n^4$ and $(n+1)^4$.
\item For each integer $n > 0$, there is a prime quadruplet between $n^5$ and $(n+1)^5$.
\item For each integer $n > 0$, there is a prime quintuplet between $n^6$ and $(n+1)^6$.
\item For each integer $n > 6$, there is a prime sextuplet between $n^7$ and $(n+1)^7$.
\end{itemize}
A further generalization is also possible:
\begin{itemize}
\item There is a prime $k$-tuplet between $n^r$ and $(n+1)^r$ for each integer $n > n_0(k, r)$, 
      where $n_0(k,r)$ is a function of $k\ge1$ and $r > 1$. 
\end{itemize}
To justify the above Legendre-type conjectures, we can assume the $k$-tuple conjecture 
plus statement (B) (sect.\,4.2) bounding the size of gaps between $k$-tuples: 
$g_k(p)<M_k \log^{k+1}p$. We can now use the following elementary argument:
Consider a fixed $r > 1$, and let $x$ be a number in the interval between $n^r$ and $(n+1)^r$. 
Then, for large $n$, the interval size $d_r=(n+1)^r-n^r\sim rn^{r-1}$ 
will be asymptotic to $rx^{(r-1)/r}$: 
because $x\sim n^r$ and $d_r \sim rn^{r-1}$ when $n \to\infty$, 
we have $n\sim x^{1/r}$ and $d_r \sim rx^{(r-1)/r}$ when $x \to\infty$. 
But any positive power of $x$ grows faster than any positive power of $\log x$ when $x \to\infty$.
So $x^{(r-1)/r}$ must grow faster than $\log^{k+1} x$. Therefore, the intervals 
$[n^r,(n+1)^r]$ --- whose sizes are about $rx^{(r-1)/r}$ --- will eventually become
much larger than the largest gaps between prime $k$-tuples containing primes $p \approx x$.
For smaller $n$, a computer check finishes the job.

However, {\bf\em this is not a proof}: we have relied on unproven assumptions.
As Hardy and Wright pointed out in 1938 (referring to the infinitude of twin primes and prime triplets),
\begin{quotation}\noindent{\small
Such conjectures, with larger sets of primes, can be multiplied, but their
proof or disproof is at present beyond the resources of mathematics. \cite[p.\ 6]{hw}
}
\end{quotation} 
Many years have passed, yet conjectures like these remain exceedingly difficult to prove.

\section{Appendix: a note on statistics of extremes}

In this appendix we use extreme value statistics to derive a simple formula expressing the 
expected maximal interval between rare random events in terms of the average interval:
$$
E(\mbox{max interval}) = a \log(T/a) + O(a)  \qquad\qquad   (1 \ll a \ll T), 
$$
where $a$ is the average interval between the rare events, 
$T$ is the total observation time or length, as applicable,
and $E(\mbox{max interval})$ stands for the {\em mathematical expectation} of the maximal interval.
The formula holds for random events occurring at exponentially distributed (real-valued) intervals, 
as well as for events occurring at geometrically distributed discrete (integer-valued) intervals. 
(For more information on extreme value distributions of random sequences see Gumbel's classical book \cite{gumbel} 
or more recent books \cite{bgt},\,\cite{galambos}.
For extreme value distributions of discrete random sequences, such as head runs in coin toss sequences, 
see also the papers of Schilling \cite{schilling} or 
Gordon, Schilling, and Waterman \cite{gsw} and further references therein.)

\subsection{Two problems about random events}
For illustration purposes, we will use two problems:

\medskip\noindent
{\bf Problem A. }
Consider a non-stop toll bridge with very light traffic. Let 
$P > 1/2$ be the probability that no car crosses the toll line during a one-second interval, and 
$q = 1 - P$ the probability to see a car at the toll line during any given second. 
Suppose we observe the bridge for a total of $T$ seconds, where $T$ is large, while $P$ is constant.

\medskip\noindent
{\bf Problem B. } Consider a biased coin with a probability of heads $P > 1/2$ (and the probability of tails $q = 1 - P$). 
We toss the coin a total of $T$ times, where $T$ is large.

\medskip\noindent
In both problems, answer the following questions about the rare events (cars or tails):

\medskip\noindent
(1) What is the {\it expected total number} of rare events 
in the observation series of length $T$?

\medskip\noindent
(2) What is the {\it expected average interval} $a$ between events (i.\,e., between cars/tails)?

\medskip\noindent
(3) What is the {\it expected maximal interval} between events, as a function of $a$?

\medskip\noindent
Notice that the first two questions are much easier than the third.
Here are the easy answers:

\medskip\noindent
(1) Because the probability of the event is $q$ at any given second/toss, 
we expect a total of $nq$ events after $n$ seconds/tosses, 
and a total of $Tq$ events at the end of the entire observation series of length $T$.

\medskip\noindent
(2) To estimate the expected average interval $a$ between events, 
we divide the total length $T$ of our observation series by the expected total number of events $Tq$. 
So a reasonable estimate\footnote{
For a small $q$, the estimate $a \approx 1/q$ is quite accurate: its error is only $O(1)$.
To prove this, we can use specific distributions of intervals between events.
Thus, if in Problem A the intervals between cars are distributed exponentially (CDF $1-P^t=1-e^{-t/a}$),
then the mean interval is $a=1/\log(1/P) = 1/q+O(1)$. 
If in Problem B the observed runs of heads are distributed geometrically (CDF $1-P^{r+1}$),
then the mean run of heads is $P/q = 1/q+O(1)$. 
}
of the expected average interval between events is $a \approx T/(Tq) = 1/q$.

\medskip\noindent
(3) Quite obviously, we can predict that the expected maximal interval is less than $T$, but not less than $a$:
$$
    a  \le  E(\mbox{max interval})  <  T. 
$$
The expected maximal interval will likely depend on both $a$ and $T$:
$$
    E(\mbox{max interval})  =  f(a,T). 
$$
It is also reasonable to expect that $f(a,T)$ should be an increasing function of both arguments, $a$ and $T$.
Can we say anything more specific about the expected maximal interval?

\subsection{An estimate of the most probable maximal interval}
In both problems A and B we will assume that $1 \ll a \ll T$ --- or, in plain English:
\begin{itemize}
\item the events are {\it rare} ($1\ll a$), and
\item our observations continue for long enough to see many events ($a \ll T$).
\end{itemize}
In Problem A, to estimate the most probable maximal interval between cars we proceed as follows:
After $n$ seconds of observations, we would have seen about $nq$ cars, 
hence about $nq$ intervals between cars. The intervals are independent of each other and real-valued. 
A known good model for the distribution of these intervals is the {\it exponential distribution}
that has the cumulative distribution function (CDF) $1-P^{\,t}$: 


with probability $P$, any given interval between cars is at least 1 second;

with probability $P^2$, any given interval is at least 2 seconds;

with probability $P^3$, any given interval is at least 3 seconds;$\ldots$

with probability $P^{\,t}$, any given interval is at least $t$ seconds.

\smallskip\noindent
Thus, after $n$ seconds of observations and about $nq$ carless intervals, we would reasonably expect 
that at least one interval is no shorter than $t$ seconds if we choose $t$ such that
$$
    P^{\,t} \times (nq) \ge 1. 
$$
Now it is easy to estimate the most probable maximal interval $t_{\max}$:
$$
    P^{\,t_{\max}}  \approx  1/(nq)
$$
$$
    (1/P)^{t_{\max}}  \approx  nq
$$
$$
    t_{\max}  \approx  \log_{1/P}(nq). 
$$
In Problem B we can estimate the longest run of heads $R_n$ after $n$ coin tosses reasoning very similarly. 
One notable difference is that now the head runs are discrete (have integer lengths). 
Accordingly, they are modeled using the {\it geometric distribution}. 
Schilling \cite{schilling} has this estimate for the longest run of heads after $n$ tosses, 
given the heads probability $P$:
$$
    R_n  \approx  \log_{1/P}(nq). 
$$
In both problems, the estimates for the most probable maximal interval (as a function of $P$ and $n$) 
have the same form $\log_{1/P}(nq)$. Therefore, it is reasonable to expect that the answers to our 
original question (3) in both problems A and B will also be the same or similar functions of 
the average interval $a$, even though the problems are modeled using different distributions of intervals. 
We will soon see that this indeed is the case.

\subsection{If random events are rare...}
If the events (cars in Problem A, or tails in Problem B) are rare, 
then $P$ is close to 1, and $q$ is close to 0. 
Using the Taylor series expansion of $\log(1/(1 - q))$, we can write:
$$
\log(1/P) = \log\left({1 \over{1-q}}\right) = q + {q^2\over 2} + {q^3\over 3} + \ldots  = q + O(q^2) 
$$
or, omitting the $O(q^2)$ terms,
\begin{eqnarray*}
\log(1/P) &\approx&  q,   \qquad\qquad\mbox{ and therefore } \\
{1 \over \log(1/P)} &\approx& {1\over q} ~\approx~ a \qquad (\mbox{moreover, we have } {1\over\log(1/P)}=a \mbox{ \,in Problem A}). 
\end{eqnarray*}
So we can transform the estimate of most probable maximal intervals, $\log_{1/P}(nq)$, like this:
$$
    \log_{1/P}(nq)  =  \log(nq) / \log(1/P)  \approx  (1/q) \log(nq) \approx  a \log(n/a). 
$$
For a long series of observations, with the total length or duration $n = T$ 
(e.\,g.~$T$ tosses of a biased coin, or $T$ seconds of observing the bridge), 
the estimate for the most probable maximal interval becomes $a \log(T/a)$.

\subsection{Expected maximal intervals}
The specific formulas for expected maximal intervals between rare events 
depend on the nature of events in the problem (whether 
the initial distribution of intervals is {\it exponential} or {\it geometric}). 
However, as $T \to\infty$, in the formulas for both cases the highest-order term turns out to be the same: 
$a \log(T/a)$, which was precisely our estimate for the most probable maximal interval.

\medskip\noindent
{\bf (A) Exponential initial distribution.}
Fisher and Tippett \cite{ft}, Gnedenko \cite{gnedenko}, Gumbel \cite{gumbel} and other authors 
showed that, for initial distributions of {\it exponential type}
(including, as a special case, the exponential distribution) 
the limiting distribution of maximal terms in a random sequence is the 
{\it double exponential distribution}~--- often called the {\it Gumbel distribution}. 
In particular, if intervals between cars in Problem A have exponential distribution 
with CDF $1-P^{\,t} = 1-e^{-t/a}$,
then the distribution of maximal intervals has these characteristics\footnote{Instead of the 
{\it scale parameter} $a$,
Gumbel \cite[p.\,157]{gumbel} uses the parameter $\alpha=1/a$. 
The {\it mode} $\mu_N$ (most probable value, also called the {\it location parameter}) in the $N$-event extreme-value distribution 
resulting from an exponential initial distribution is equal to the 
{\it characteristic extreme} $a \log N$ \cite[p.\,114]{gumbel}. 
The shape of the $N$-event extreme-value distribution approaches that of the limiting distribution as $N\to\infty$.
}:
\begin{eqnarray*}
\mbox{$N$-event CDF:}&\hphantom{=}& (1-e^{-t/a})^N = (1-\textstyle{1\over N}e^{-(t-\mu_N)/a})^N 
\quad\mbox{(distribution for $N\approx Tq$ events), } \\
\mbox{Limiting CDF:} &\hphantom{=}& \exp(-e^{-(t-\mu)/a})\qquad\mbox{(Gumbel distribution) \cite[{\rm p.\,157}]{gumbel},}\\
\mbox{Scale} &=&  a   ~=~ 1/\log(1/P) \quad\mbox{ (equal to the expected average interval)},\\
\mbox{Mode}  &=&  \mu  ~=~ \mu_N ~=~ a \log N ~\approx~a \log(T/a) ~\approx~\log_{1/P}(Tq) , \\
\mbox{Median}  &=&  \mu-a\log\log2 ~\approx~ a \log(T/a) + 0.3665a, \\
\mbox{Mean}  &=&  \mu + \gamma a ~\approx~ a \log(T/a) + 0.5772a,
\end{eqnarray*}
where $\gamma = 0.5772...$ is the Euler-Mascheroni constant. 
The mean value of observed maximal intervals in Problem A will converge almost surely 
to the mean $\mu + a\gamma$ of the Gumbel distribution, therefore:
$$
E(\mbox{max interval}) ~\approx~ \log_{1/P}(Tq) + \gamma a ~\approx~ a\log(T/a) + \gamma a   ~=~  a \log(T/a) + O(a). 
$$

\medskip\noindent
{\it Historical notes: }  
In 1928 Fisher and Tippett \cite{ft} described three types of limiting extreme-value distributions 
and showed that the double exponential (Gumbel) distribution is the limiting extreme-value distribution 
for a certain wide class of random sequences. They also computed, among other parameters, 
the {\it mean-to-mode distance} in the double exponential distribution \cite[p.\,186]{ft}; 
it is this result that allows one to conclude that the mean is $\mu + a\gamma$ if the mode is $\mu$.
Gnedenko (1943) \cite{gnedenko} rigorously proved the necessary and sufficient conditions 
for an initial distribution to be in the domain of attraction of a given type of limiting distribution. 

\medskip\noindent
{\bf (B) Geometric initial distribution. }
Surprisingly, in this case the limiting extreme-value distribution does not exist 
\cite[p.\,203]{schilling},\,\cite[p.\,280]{gsw}. 
For the longest run of heads $R_n$ in a series of $n$ tosses of a biased coin,
with the probability of heads $P$, we have 
$$
E(R_n)  =  \log_{1/P}(nq) + \frac{\gamma}{\log(1/P)} - {1\over2} + \mbox{smaller terms}  
\qquad\mbox{\cite[{\rm p.\,202}]{schilling}},
$$
where the first term is the same as in Problem A (up to a substitution $n = T$). 
The sum of the other terms is $O(a)$ when $P$ is close to 1; 
so, again, we have
$$
    E(R_n)  =  a \log(n/a) + O(a). 
$$
 
\subsection{Standard deviation of extremes}
As above, the specific formula for standard deviation (SD) in distributions of 
maximal intervals between events depends on the nature of the problem 
(whether the initial distribution of intervals is exponential or geometric). 
Still, in both cases SD $\approx \pi a/\sqrt{6} = O(a)$.

\medskip\noindent
{\bf (A) Exponential initial distribution.}
Here the limiting distribution of maximal intervals is the Gumbel distribution with the scale $a = 1/\log(1/P)$, 
therefore the SD of maximal intervals must be very close to the SD of the Gumbel distribution:
$$
\mbox{SD}(\mbox{max interval}) \approx {\pi a \over \sqrt{6}} =O(a)
 \qquad\mbox{\cite[{\rm p.\,116,\,174}]{gumbel}}.
$$

\medskip\noindent
{\bf (B) Geometric initial distribution.}
For the longest run of heads $R_n$ in a series of $n$ tosses of a biased coin, the variance is
$$
\mbox{Var\,}R_n = {\pi^2\over 6\log^2(1/P)} + {1\over12} + \mbox{smaller terms} 
\qquad\mbox{\cite[{\rm p.\,202}]{schilling}},
$$
where the first term is $O(a^2)$, 
while the sum of the other terms is much smaller than the first term. 
(Again, recall that for average intervals $a$ between rare events~--- in this case, 
between {\it tails}~--- we have $a \approx 1/\log(1/P)$.)
Therefore, the standard deviation is 
$$
\mbox{SD}\,\,R_n  = \sqrt{\mbox{Var\,}R_n} = {\pi \over \sqrt{6} \log(1/P)} + \mbox{a small term} \approx {\pi a \over\sqrt{6}}=O(a).
$$

\subsection{A shortcut to the answer}
There is a simple way to ``guesstimate'' the answer $a\log(T/a) + O(a)$. 
If $a$ is the average interval between events, 
then the most probable maximal interval is about $a\log(T/a)$ (sect.~8.3). 
We can now simply use the fact that the width of the extreme value distribution is $O(a)$. 
(Imagine what happens if the rare event's probability $q$ is reduced by 50\%.
This change in $q$ would have about the same effect as if every interval became twice as large: 
then average and maximal intervals would also become twice as large, 
and the extreme value distribution would be twice as wide. 
This immediately implies that the extreme value distribution is $O(a)$ wide.) 
But then the true value of the expected maximal interval cannot be any farther 
than $O(a)$ from our estimate $a\log(T/a)$; so the expected maximal interval is $a\log(T/a) + O(a)$.

\subsection{Summary}
We have considered maximal intervals between random events in two common situations:
\begin{itemize}
\item rare events occurring at {\it exponentially distributed} intervals (Problem A);
\item discrete rare events at {\it geometrically distributed} intervals (Problem B).
\end{itemize}
These two situations are somewhat different:
in the former case maximal intervals have a limiting distribution (the Gumbel distribution), while
in the latter case no limiting distribution exists (here the Gumbel distribution is simply a decent approximation).
Nevertheless, in both cases the expected maximal interval between events is
$$
E(\mbox{max interval}) = a \log(T/a) + \gamma a + \mbox{lower-order terms } =  a \log(T/a) + O(a),    
$$
where $a$ is the average interval between events, $T$ is the total observation time or length,
and the lower-order terms depend on the initial distribution. 

\medskip\noindent
As we have seen in Sections 4--6, a remarkably similar {\it heuristic formula} $a\log(x/a)-ba$, 
with an empirical term $-ba$ replacing the ``theoretical'' $\gamma a$, satisfactorily describes the following:
\begin{itemize}
\item record gaps between primes below $x$ ($a = \log x$, $b \approx 3$; A005250)
\item record gaps between twin primes below $x$ ($a = 0.75739 \log^2 x$, $b\approx1$; A113274) and, more generally,
\item record gaps between prime $k$-tuples ($a = C_k \log^k x$, $b \approx 2/k$, 
      where $C_k$ is reciprocal to the Hardy-Littlewood constant for the particular $k$-tuple).
\end{itemize}
%

\section{Acknowledgements}
The author is grateful to the anonymous referee and to
all authors, contributors, and editors of the websites 
{\it OEIS.org}, {\it PrimePuzzles.net} and {\it FermatQuotient.com}.
Many thanks also to Prof.~Marek Wolf for his interest in the initial version of this paper, 
followed by an email exchange that undoubtedly helped make this paper better.

\bigskip
\hrule
\bigskip

\noindent 2010 {\it Mathematics Subject Classification}:
Primary 11N05; Secondary 60G70.

\noindent \emph{Keywords: } 
distribution of primes, prime $k$-tuple, Hardy-Littlewood conjecture,
extreme value statistics, Gumbel distribution,
prime gap, Cram\'er conjecture, 
prime constellation, twin prime conjecture, prime quadruplet, prime sextuplet.

\bigskip
\hrule
\bigskip

\noindent (Concerned with OEIS sequences A005250, A091592, A113274, A113404, A192870, A200503, 
A201596, A201598, A201051, A201062, A201073, A201251, A202281, A202361.) 

\bigskip
\hrule
\bigskip

\vspace*{+.1in}
\noindent


\noindent
\vskip .1in


\begin{thebibliography}{20}


\bibitem{bgt}
J.~Beirlant, Y.~Goegebeur, J.~Segers and J.~Teugels, {\it Statistics of Extremes: Theory and Applications}, Wiley, 2004.

\bibitem{berman}
S.~M.~Berman, Limiting distributions of the maximum term in sequences of dependent random variables.
{\it Annals of Mathematical Statistics}, 33 (1962), 894--908.

\bibitem{cohen}
H.~Cohen, High precision computation of Hardy-Littlewood constants, {preprint. \hfill} \linebreak
{\small\tt http://www.math.u-bordeaux1.fr/\~{}cohen/hardylw.dvi} (2012).


\bibitem{cram}
H.~Cram\'er, On the order of magnitude of the difference between consecutive prime numbers. 
{\it Acta Arith.} 2 (1937), 23--46. 

\bibitem{fischer} 
R.~Fischer, Maximale L\"ucken (Intervallen) von Primzahlenzwillingen, web page (in German) 
{\small\tt http://www.fermatquotient.com/PrimLuecken/ZwillingsRekordLuecken} (2008).

\bibitem{fischer2} 
R.~Fischer, Maximale Intervalle von Primzahlenpaaren, web page (in German){. \hfill} \linebreak
Maximal gaps between twin primes are $G_2(p) \approx (1.32032)^{-1}(\log p-(2/3) \log \log p)^3${. \hfill} \linebreak
Maximal gaps between prime triplets are $G_3(p) \approx (2.8582)^{-1}(\log p-(3/4) \log \log p)^4${. \hfill} \linebreak
{\small\tt http://www.fermatquotient.com/PrimLuecken/Max\_Intervalle} (2006).

\bibitem{ft}
R.~A.~Fisher and L.\,H.\,C.~Tippett, Limiting forms of the frequency distribution of 
the largest and smallest member of a sample, {\it Proc. Camb. Phil. Soc.}, 24 (1928), 180--190.

\bibitem{forbes}
A.~D.~Forbes, Prime $k$-tuplets. Section 21: List of all possible patterns of prime $k$-tuplets. \linebreak
The Hardy-Littlewood constants pertaining to the distribution of prime $k$-tuplets{. \hfill} \linebreak
{\small\tt http://anthony.d.forbes.googlepages.com/ktuplets.htm} (2012).

\bibitem{galambos}
J.~Galambos, {\it The Asymptotic Theory of Extreme Order Statistics}, 
Krieger, 1987. 

\bibitem{gnedenko}
B.~V. Gnedenko, Sur la distribution limite du terme maximum d'une s\'erie al\'eatoire. {\it Ann. Math.}, 44 (1943), 423--453.
(English translation: On the limiting distribution of the maximum term in a random series. 
{\it Breakthroughs in Statistics}, Volume 1: Foundations and Basic Theory. Springer, New York, 1993, pp.\,185--225.)





\bibitem{gsw}
L.~Gordon, M.~F.~Schilling, and M.~S.~Waterman, An extreme value theory for long head runs. 
{\it Probability Theory and Related Fields}, 72 (1986), 279--297.

\bibitem{granville}
A. Granville. Harald Cram\'er and the distribution of prime numbers. 
{\it Scand. Act. J.}, 1 (1995), 12--28.

\bibitem{gumbel}
E.~J.~Gumbel, {\it Statistics of Extremes}, Columbia University Press, 1958. Dover, 2004.

\bibitem{hl}
G.~H.~Hardy and J.~E.~Littlewood, Some Problems of `Partitio Numerorum.' III. 
On the Expression of a Number as a Sum of Primes. {\it Acta Math.} 44 (1922), 1--70.

\bibitem{hw} G. H. Hardy and E. M. Wright, 
{\it An Introduction to the Theory of Numbers}, 6th ed. Oxford University Press, 2008. 

\bibitem{hays}
W.~L.~Hays, {\it Statistics}, Harcourt Brace College Publishers, 1994.

\bibitem{kpII}
P.~F.~Kelly and T.~Pilling,
Implications of a new characterisation of the distribution of twin primes.
{\small\tt http://arxiv.org/abs/math/0104205} (2001).

\bibitem{kpIII}
P.~F.~Kelly and T.~Pilling,
Physically inspired analysis of prime number constellations.
{\small\tt http://arxiv.org/abs/hep-th/0108241} (2001).

\bibitem{maier}
H.~Maier, Primes in short intervals, {\it Michigan Math. J.}, 32 (1985), 221--225.

\bibitem{mp}
H.~Maier and C.~Pomerance, Unusually large gaps between consecutive primes. 
{\it Transactions of the AMS}, 322 (1990), 201--237. 

\bibitem{easyfit}
MathWave Technologies, {\it EasyFit --- Distribution Fitting Software}, the company web site at 
{\small\tt http://www.mathwave.com/easyfit-distribution-fitting.html} (2012).

\bibitem{exceltrendline}
Microsoft Corporation, {\it Excel: Add a Trendline to a Chart}, the company web site {at \hfill} \linebreak
{\small http://office.microsoft.com/en-us/excel-help/add-a-trendline-to-a-chart-HP005198462.aspx} (2012).

\bibitem{nicely}
T.~R.~Nicely, List of prime gaps. {\small\tt http://www.trnicely.net/gaps/gaplist.html} (2012).

\bibitem{cvc}
J.~Pintz, Cram\'er vs Cram\'er: On Cram\'er's probabilistic model of primes.
{\it Functiones et Approximatio}, XXXVII.2 (2007), 361--376.


\bibitem{riesel}
H.~Riesel, {\it Prime Numbers and Computer Methods for Factorization}, 
Birkh\"auser, 1994. 

\bibitem{rr} 
L. Rodriguez and C. Rivera, Conjecture 66. Gaps between consecutive twin prime pairs.
{\small\tt http://www.primepuzzles.net/conjectures/conj\_066.htm} (2009).

\bibitem{schilling}
M. F. Schilling, The longest run of heads. {\it The College Math. J.}, 21 (1990), 196--207. 

\bibitem{shanks}
D.~Shanks, On maximal gaps between successive primes. {\it Math. Comput.} 18 (1964), 646--651. 

\bibitem{oeis}
N.\,J.\,A.\,Sloane, {\it On-Line Encyclopedia of Integer Sequences}, 
{\small\tt http://oeis.org} (2012).

\bibitem{watson}
G.~S.~Watson, Extreme values in samples from $m$-dependent stationary stochastic processes.
{\it Annals of Mathematical Statistics}, 25 (1954), 798--800.

\bibitem{wolf} 
M.~Wolf, Some heuristics on the gaps between consecutive primes, arXiv preprint.
{\small\tt http://arxiv.org/abs/1102.0481} (2011)

\bibitem{wolf2} 
M.~Wolf, Maximal gaps between twin primes $G_2(x)$ can be expressed in terms of $\pi_2(x)$.
E-mail communication (2013).
\end{thebibliography}
\end{document}